\documentclass[12pt]{article}
\usepackage{amsfonts,amssymb,epsfig}
\usepackage{amsmath}

\date{}
\begin{document}
\newtheorem{df}{Definition}
\newtheorem{thm}{Theorem}

\newtheorem{lm}{Lemma}
\newtheorem{pr}{Proposition}
\newtheorem{co}{Corollary}
\newtheorem{re}{Remark}
\newtheorem{note}{Note}
\newtheorem{claim}{Claim}
\newtheorem{problem}{Problem}

\def\R{{\mathbb R}}

\def\E{\mathbb{E}}
\def\calF{{\cal F}}
\def\N{\mathbb{N}}
\def\calN{{\cal N}}
\def\calH{{\cal H}}
\def\n{\nu}
\def\a{\alpha}
\def\d{\delta}
\def\t{\theta}
\def\e{\varepsilon}
\def\t{\theta}
\def\g{\gamma}
\def\G{\Gamma}
\def\b{\beta}
\def\pf{ \noindent {\bf Proof: \  }}

\newcommand{\qed}{\hfill\vrule height6pt
width6pt depth0pt}
\def\endpf{\qed \medskip} \def\colon{{:}\;}
\setcounter{footnote}{0}

\renewcommand{\qed}{\hfill\vrule height6pt  width6pt depth0pt}

\title{Polylog dimensional subspaces of $\ell_\infty^N$\thanks {2010 AMS subject classification: 46B07, 46B20.\
Key words: the dichotomy problem, $\ell_\infty^n$,  cotype.
This research was completed in Fall 2017 while both authors where members of the Geometric Functional Analysis and Application program at MSRI, supported by the National Science Foundation under Grant No. 1440140}
}

\author{Gideon Schechtman\thanks{Supported in part by the Israel Science
Foundation
 } and  Nicole Tomczak--Jaegermann} \maketitle

\begin{abstract}
We show that a subspace of of $\ell_\infty^N$ of dimension $n>(\log N\log \log N)^2$ contains $2$-isomorphic copies of $\ell_\infty^k$ where $k$ tends to infinity with $n/(\log N\log \log N)^2$. More precisely, for every $\eta>0$, we show that any subspace of $\ell_\infty^N$ of dimension $n$ contains a subspace of dimension \hfil\break
$m=c(\eta)\sqrt{n}/(\log N\log \log N)$ of distance at most $1+\eta$ from $\ell_\infty^m$.
\end{abstract}

\section{Introduction}

The dichotomy problem of Pisier asks whether a Banach space $X$  either contains, for  every $n$, a subspace $K$-isomorphic to $\ell_\infty^n$, for some (equivalently all) $K>1$, or, for every $n$, every $n$-dimensional subspace of $X$ $2$-embeds in $\ell_\infty^N$ only if $N$ is exponetial in $n$. This is equivalent to the question of whether for some (equivalently all) absolute $K>1$ and any sequence $n_N\le N$ with $n_N/\log N\to \infty$ when $N\to\infty$, every subspace of $\ell_\infty^N$ of dimension $n_N$ contains a subspace of dimension $m_N$ $K$-isomorphic to $\ell_\infty^{m_N}$ where $m_N\to\infty$ when $N\to\infty$.

We remark in passing that the equivalence between the two versions of the problem (``some $K>1$" versus ``all $K>1$") is due to the fact proved by R.C. James that, for all $1<\kappa<K<\infty$, a space which is $K$ isomorphic to $\ell_\infty^n$ contains a subspace $\kappa$ isomorphic to $\ell_\infty^m$ where $m\to\infty$ as $n\to\infty$.

As is exposed in \cite{p}, Maurey proved that if $X^*$ has non-trivial type (Equivalently does not contain uniformly isomorphic copies of $\ell_1^n$-s. This is a condition stronger than $X$ has non-trivial cotype; equivalently, does not contain uniformly isomorphic copies of $\ell_\infty^n$-s), then we get the required conclusion: For every $n$ every $n$-dimensional subspace of $X$ $2$-embeds in $\ell_\infty^N$ only if $N$ is exponetial in $n$.

Another partial result was obtained by Bourgain in \cite{b2} where he showed in particular that the conclusion holds if $n_N>(\log N)^4$.

Here we show some improvement over this result of Bourgain: The conclusion holds if  $n_N/(\log N\log\log N)^2$ tends to $\infty$.

\begin{thm}\label{thm:main}
Let $n$, and $N$ be integers
such that $n>(\log N\log \log N)^2$. Then, for some absolute constant $c>0$ and for every $0<\eta<1$, any subspace of $\ell_\infty^N$ of dimension $n$ contains a subspace of dimension \hfil\break
$m=c\eta^2\sqrt{n}/(\log N\log \log N)$ of distance at most $1+\eta$ from $\ell_\infty^m$.
\end{thm}

Note that we get some specific estimates for the dimension of the contained subspace $(1+\eta)$-isomorphic to an $\ell_\infty$ space of its dimension. Although we are interested in small $n$-s, the result gives some estimate in the whole range. This is also the case in Bourgain's result: He proved that if $n\ge N^\delta$ than any subspace of $\ell_\infty^N$ of dimension $n$ contains a subspace $(1+\eta)$-isomorphic to an $\ell_\infty$ of dimension $m\ge c\eta^5\delta^2\sqrt{n}/\log (1/\delta)$.  Comparing the two, our result gives better estimates for $m$ when $n\lesssim e^{c(\eta)\sqrt{\log N}}$ and worse when $n$ is larger.
Recall also that for $n$ proportional to $N$, Figiel and Johnson \cite{fj} proved earlier that $m$ can be taken of order $\sqrt{N}$ (and no better). This is not recovered by our result.

The general idea of the proof of Theorem \ref{thm:main} is the same as in \cite{b2} but the technical details are somewhat different. At the end of this note we also speculate that, up to the $(\log\log N)^2$ factor, our result may be best possible.

Our result was essentially achieved a long time ago, circa 1990. Since several people showed interest in it lately we decided to write it up with the hope that more modern methods (and younger minds) may be able to improve it farther.

 \section{Proofs}\label{section:proofs}

 The main technical tool in the proof of Theorem \ref{thm:main} is the following proposition

\begin{pr}\label{pr:main}
Let $n$, and $N$ be integers
such that $n>(\log N)^{3/2} \log \log N$.
Let $[a_i(j)]$ be an $n \times N$ matrix with $a_i(j)\ge 0$ for
$i= 1, \dots, n$ and $j = 1, \dots, N$.
 Assume that
\[
\sum_{i=1}^n  a_i(j)^2 \le 1 \ \ \mbox {for  }\ \   j = 1, \dots, N
\]
and
\[
\sum_{i=1}^n  a_i(j) \le 3\sqrt{log N} \ \ \mbox{ for  }\ \  j = 1, \dots, N.
\]
Moreover, assume that, for some $\gamma>0$,
for every $i= 1, \dots, n$ there exists  $1 \le j \le N$
such that $a_i(j) \ge \gamma$.
Denote by $a_i$ the $i$-th row of the matrix. Then, for some positive constants, $c(\gamma), K(\gamma)$ depending only on $\gamma$ and for every $0<\eta<1$,  there are disjoint subsets $\sigma_1,\dots,\sigma_m$ of $\{1,\dots,n\}$ with $m \ge c(\gamma)\eta^2 n/(\log N)^{3/2}
 \log\log N$,
 Such that
 \[
 \|\sum_{r=1}^m\sum_{i\in\sigma_r} a_i\|_\infty/\min_{1\le r\le m}\|\sum_{i\in\sigma_r} a_i\|_\infty\le (1+K(\gamma)\eta).
 \]
 \end{pr}

 We first show how to deduce Theorem \ref{thm:main} from the proposition above.

 \noindent {\bf Proof of Theorem \ref{thm:main}:} Let $X$ be an $n$ dimensional subspace of $\ell_\infty^N$. The $\pi_2$ norm of the identity on $X$ is equal to $\sqrt{n}$ (\cite{gg},\cite{s}) and by the main theorem of \cite{t-j} (see \cite{t-j1} for the constant $\sqrt 2$) this quantity can be computed, up to constant $\sqrt 2$ on $n$ vectors. This means that there are $n$ vectors $a_i=(a_i(1),\dots,a_i(N))$, $i=1,\dots,n$, in $X$ satisfying
 \[
 \sum_{i=1}^n a_i(j)^2\le 1, \ \ \mbox{for all}\ \  j=1,\dots,N
 \]
 and
 \[
 \sum_{i=1}^n\|a_i\|_\infty^2\ge n/2.
 \]
 The first condition implies in particular that $\|a_i\|_\infty^2\le 1$ for each $i$ so necessarily for a subset $\sigma^\prime$ of $\{1,\dots,n\}$ of cardinality at least $n/4$, $\|a_i\|_\infty\ge 1/2$ for all $i\in \sigma^\prime$. The existence of a subset $\sigma^\prime$ of $\{1,\dots,n\}$ of cardinality at least $n/4$ satisfying the two conditions
 \begin{equation}\label{eq:2summ}
 \sum_{i\in \sigma^\prime}a_i(j)^2\le 1, \ \ \mbox{for all}\ \  j=1,\dots,N, \ \ \mbox{and}\ \ \|a_i\|_\infty\ge 1/2 \ \ \mbox{for all}\ \ i\in\sigma^\prime
 \end{equation}
 is all that we shall use from now on. In Remark \ref{re:alternate 2summ} below we'll show another way to obtain this.

 Next we would like to choose a subset $\sigma$ of $\sigma^\prime$ of cardinality of order $\sqrt{n\log N}$ such that the matrix $[|a_i(j)|]$, $i\in\sigma$, $j=1,\dots,N$, will satisfy the assumptions of Proposition \ref{pr:main}.
 So let $\xi_i$, $i\in\sigma^\prime$, be independent $\{0,1\}$ valued random variables with ${\rm Prob}(\xi_i=1)=\sqrt{(\log N)/n}$. Since for all $j$ $\sum_{u\in\sigma^\prime}|a_i(j)|\le \sqrt n$, $\E\sum_{u\in\sigma^\prime}|a_i(j)|\xi_i\le \sqrt{\log N}$. By the most basic concentration inequality, using the fact that $\sum_{i\in\sigma^\prime}a_i(j)^2\le 1$, for all $j$,
 \begin{multline*}
 {\rm Prob}(\sum_{i\in\sigma^\prime}|a_i(j)|\xi_i>3\sqrt{\log N})\\
 \le {\rm Prob}(\sum_{i\in\sigma^\prime}|a_i(j)|(\xi_i-\E\xi_i)>2\sqrt{\log N})\le e^{-2\log N}=1/N^2.
 \end{multline*}
 It follows that with probability larger than $1-1/N$
 \[
 \sum_{i\in\sigma^\prime}|a_i(j)|\xi_i\le3\sqrt{\log N}
 \]
 for all $j$. Since by a similar argument also $\sum_{i\in\sigma^\prime}\xi_i\ge \frac{\sqrt{n\log N}}{16}$ with probability tending to $1$ when $N\to\infty$ we get a subset $\sigma$ of  cardinality $n^\prime\ge \frac{\sqrt{n\log N}}{16}$ satisfying
 \[
 \sum_{i\in\sigma}|a_i(j)|\le3\sqrt{\log N}\ \ \mbox{for all}\ \ j=1,\dots,N.
 \]
 Note that the condition $n\ge 256(\log N\log\log N)^2$ implies that \hfil\break
 $n^\prime\ge (\log N)^{3/2}\log\log N$. It follows that the matrix $[|a_i(j)|]$, $i\in \sigma^\prime, j=1,\dots,N$ satisfies the conditions of Proposition \ref{pr:main} with $n^\prime$ replacing $n$ and $\gamma=1/2$. We thus get that, for some absolute positive constants $c,K$, there are disjoint subsets $\sigma_1,\dots,\sigma_m$ of $\{1,\dots,n\}$ with
 \[
 m \ge 16c\eta^2 n^\prime/(\log N)^{3/2}\log\log N \ge c\eta^2\sqrt{n}/\log N\log\log N,
 \]
 such that
 \[
 \|\sum_{r=1}^m\sum_{i\in\sigma_r} |a_i|\|_\infty/\min_{1\le r\le m}\|\sum_{i\in\sigma_r} |a_i|\|_\infty\le (1+K\eta).
 \]
 Rescaling, we may assume that $\min_{1\le r\le m}\|\sum_{i\in\sigma_r} |a_i|\|_\infty=1$. Let $j_r$ denote the label of (one of) the largest coordinates of $\sum_{i\in\sigma_r} |a_i|$. Assume as we may that $\eta<1/K$. Then no two $r$'s can share the same $j_r$. Changing the labelling we can also assume $j_r=r$.

 Put $x_r= \sum_{i\in\sigma_r}sign(a_i(r))a_i$. Then for all $r$, $\|x_r\|_\infty\ge 1$ and for all $j=1,\dots,N$,
 \begin{equation}\label{eq:bnd}
 \sum_{r=1}^m |x_r(j)|\le 1+K\eta.
 \end{equation}
 So the sequence $x_r$, $r=1,\dots,m$, is $(1+K\eta)$-dominated by the $\ell_\infty^m$ basis; i.e,
 \[\|\sum_{r=1}^m\a_r x_r\|\le (1+K\eta)\max_{1\le r\le m}|\a_r|\ \ \mbox{for all}\ \ \{a_r\}_{r=1}^m.
 \]
 The lower estimate is achieved similarly: Assume $\max_{1\le r\le m}|\a_r|=|\a_{r_0}|$ and note that
 \[
 \|\sum_{r=1,r\not=r_0}^m\sum_{i\in\sigma_r}|a_i(r_0)|\|\le K\eta.
 \]
 Then,
 \begin{multline*}
 \|\sum_{r=1}^m\a_r x_r\|\ge |\sum_{r=1}^m\a_r x_r(r_0)|\\
 \ge |\a_{r_0}|\sum_{i\in\sigma_{r_0}}|a_i(r_0)|-\sum_{r=1,r\not=r_0}^m|\a_r|\sum_{i\in\sigma_r}|a_i(r_0)|
 \\
 \ge( (1-K\eta)\max_{1\le r\le m}|\a_r|.
 \end{multline*}
 We have thus found a subspace of $x$ of dimension $m\ge c\eta\sqrt{n}/(\log N\log\log N)$ whose distance to $\ell_\infty^m$ is at most $(1+K\eta)/(1-K\eta)$. Changing the last quantity to $1+\eta$, paying by changing $c$ to another absolute constant, is standard.
 \endpf

 In the proof of Proposition \ref{pr:main} we shall use the following Lemma which follows immediately from Lemma 2 in \cite{b3} (but, following the proof of that lemma from \cite{b3}, is a bit easier to conclude).

 \begin{lm}\label{lm:bourgain}Let $\xi_i$, $i\in\{1,\dots,n\}$, be independent $\{0,1\}$ valued random variables with ${\rm Prob}(\xi_i=1)=\d$. Then for all $q\ge 1$,
\[
(\E
 (\sum_{i=1}^n\xi_i)^q)^{1/q}\le C(\d n+q).
 \]
 C is a universal constant.
 \end{lm}

 We now pass to the

 \noindent {\bf Proof of Proposition \ref{pr:main}:} We shall assume as we may that $\eta<\gamma$. We first deal with the small $a_i(j)$-s.
 Fix $\varepsilon >0$ to be defined later. Let
\[
b_i(j) = \begin{cases} a_i(j) &\mbox{if } a_i(j)\le\e \\
0 & \mbox{otherwise} \end{cases}
\]
We will show that for any $\delta >0$, and for
a random subset $\sigma \subset
\{1, \dots, n\}$ of cardinality $|\sigma| \sim \delta n$
\begin{equation}\label{eq:smallcoord}
\sum_{i\in\sigma}  b_i(j) \le C(\delta \sqrt{\log N}  + \e \log N)
\ \ \mbox{for} \ \  j = 1, \dots, N,
\end{equation}
where $C$ is an absolute constant.

Indeed,
set $p = \log N$. Fix $\delta >0$ and let
$\xi_i$ denote selectors with mean $\delta$ as in Lemma \ref{lm:bourgain}.
By Chebyshev inequality, (\ref{eq:smallcoord}) follows from the estimate
\begin{equation}\label{eq:integral}
\sup_{j}
\left( \E (  \sum_{i=1}^n \xi_i(\omega) b_i(j) )^p \right)^{1/p}
\le  C(\delta \sqrt{\log N}  + \e \log N).
\end{equation}
Indeed,
\begin{multline*}
\left( \E\sum_{j=1}^N (  \sum_{i=1}^n \xi_i(\omega) b_i(j) )^p \right)^{1/p}\\
\le N^{1/\log N}
\sup_{j}\left( \E (  \sum_{i=1}^n \xi_i(\omega) b_i(j) )^p \right)^{1/p}\le  eC(\delta \sqrt{\log N}  + \e \log N).
\end{multline*}
Now apply Chebyshev's inequality.

Fix $1 \le j \le N$ and denote  $(b_i(j))_i \in \R^{l}$ by $b$.
Considering the level sets
of $b$ we may assume
without loss of generality
that $b$ is of the form
\[
b = \sum_{m=2 \log (1/\varepsilon)}^{\infty}
2^{-m/2} \chi_{D_m},
\]
($\log$ is $\log_2$) where the sets $D_m \subset \{1, \dots, n\}$ are mutually
disjoint and   $ \chi_{D_m}$ denotes the characteristic
function of the set $D_m$, for
$m = 2 \log (1/\varepsilon), \dots $.
Thus,
\begin{align}\label{eq:intestim}
&\left( \E (  \sum_{m=2 \log(1/\e)}^n \xi_j(\omega) b_i(j))^p \right)^{1/p}\nonumber
 \\
&\phantom{aaaaaaaaaaaaaa}\le  \sum_{m=2 \log(1/\e)}
^\infty 2^{-m/2}
\left( \E (  \sum_{j \in D_m} \xi_j(\omega))^p \right)^{1/p}\nonumber
\\
&\phantom{aaaaaaaaaaaaaa}\le C \sum_{m=2 \log(1/\e)}
^\infty 2^{-m/2}
\left(\delta |D_m|  + p)\right)\ \ \mbox{by Lemma \ref{lm:bourgain}}\nonumber
\\
&\phantom{aaaaaaaaaaaaaa}\le C  \delta \sum_{m=2 \log(1/\e)}
^\infty 2^{-m/2}|D_m|
+ Cp \sum_{m=2 \log(1/\e)}^\infty
2^{-m/2}
\end{align}

To estimate the first term in (\ref{eq:intestim}) note that
\[
 \sum_{m=2 \log(1/\e)}^\infty
2^{-m/2}|D_m|
= \| b \|_1 \le 3\sqrt{\log N}.
\]
The second term is clearly smaller than an absolute constant times $\e p$.

Combining the latter two estimates with (\ref{eq:intestim})
we get (\ref{eq:integral}) and hence also (\ref{eq:smallcoord}).

To deal with the large coordinates, set, for $j=1,\dots,N$,
\[
A_j=\{1\le i\le n;\ a_i(j)\ge\e\}.
\]
Since $\sum_{i=1}^n a_i(j)\le 3\sqrt{\log N}$,
\begin{equation}\label{eq:setlc}
|A_j|\le\sqrt{\log N}/\e\ \ \mbox{for}\ \ j=1,\dots,N.
\end{equation}
An argument similar to the one that proved (\ref{eq:integral}) also shows that a random set $\sigma\subset\{1,\dots,n\}$ of cardinality $|\sigma|\sim\delta n$ satisfies
\begin{equation}\label{eq:ransetlc}
|\sigma\cap A_j|\le C(\delta\sqrt{\log N}/\e +\log N),.
\end{equation}
Indeed, this follows easily by applying the following inequality with $p=\log N$,
\[
\left(\E\left(\sum_{i=1}^{\sqrt{\log N}/\e}\xi_i\right)^p\right)^{1/p}\\
\le C(\delta\sqrt{\log N}/\e+\log N).
\]
Moreover, Chebyshev's inequality implies that we can find a set $\sigma\subset\{1,\dots,n\}$ of cardinality at least $\frac12\delta n$ which satisfies (\ref{eq:smallcoord}) and (\ref{eq:ransetlc}) simultaneously (say, with the same absolute constant $C$).

Choose now $\delta=2\eta/\sqrt{\log N}$ and $\e=\eta/\log N$. Then we get a set $\sigma\subset\{1,\dots,n\}$ of cardinality at least $\eta n/\sqrt{\log N}$. such that
\begin{equation}\label{eq:small}
\sum_{i\in\sigma}b_i(j)\le 3C\eta  \ \ \mbox{for}\ \ j=1,\dots,N,
\end{equation}
and
\begin{equation}\label{eq:random}
|\sigma\cap A_j|\le 3C\log N  \ \ \mbox{for}\ \ j=1,\dots,N.
\end{equation}
define $j_1\in\{1,\dots,N\}$ and $s_1$ by
\[
s_1=\sum_{i\in \sigma\cap A_{j_1}}a_i(j_1)=\max_{j}\sum_{i\in \sigma\cap A_j}a_i(j).
\]
For $r>1$ define $S_{r-1}=\sigma\setminus(A_{j_1}\cup\dots\cup A_{j_{r-1}})$ and $j_r$ and $s_r$ by
\[
s_r=\sum_{i\in S_{r-1}\cap A_{j_r}}a_i(j_r)=\max_{j}\sum_{i\in S_{r-1}\cap A_{j}}a_i(j).
\]
By rearranging the columns we may assume $j_r=r$ for all $r$. Now, (\ref{eq:random}) implies that $| S_r|\ge |\sigma|-3Cr\log N$ so $ S_{r}$ is not empty for $1\le r\le \frac{\eta n}{3C(\log N)^{3/2}}$. Also,
\[
\gamma\le s_r\le 3\log N\ \ \mbox{for}\ \ 1\le r\le \frac{\eta n}{3C(\log N)^{3/2}}.
\]
The sequence $s_r$ is non-increasing, divide it into $(\log((3\log N)/\gamma))/\log(1+\eta)$ intervals such that in each interval $\max s_r/\min s_r$ is at most $1+\eta$.  There is an interval $R$ with $|R|\ge \frac{(\log(1+\eta))\eta n}{3C(\log N)^{3/2}\log((3\log N)/\gamma)}\ge
\frac{\eta^2 n}{6C(\log N)^{3/2}\log((3\log N)/\gamma)}$ such that
\[
\max_{r\in R} s_r/\min_{r\in R} s_r\le 1+\eta.
\]
Put $\sigma_r=\tilde S_{r-1}\cap A_r$.
Since $\min_{r\in R} s_r \ge \gamma >\eta$  we are done in view of (\ref{eq:small}) and the fact that $s_r\ge \sum_{\tilde S_{r-1}\cap A_s}a_i(s)$ for $r<s$.
\endpf

\section{Remarks}
\begin{re}\label{re:alternate 2summ}
Here is an alternative way to get (\ref{eq:2summ}):
\end{re}
Let $X$ be an $n$-dimensional normed space which, without loss of generality we assume is in John's position, i.e., the maximal volume ellipsoid inscribed in the unit ball of $X$ is the canonical sphere $S^{n-1}$. A weak form of the Dvoretzky--Rogers lemma asserts that there are orthonormal vectors $x_1,\dots,x_n$  such that $\|x_i\|_X\ge c$ for some universal positive constant $c$. This is proved by a simple volume argument, see for example Theorem 3.4 in \cite{ms}. (There it is shown that there are $[n/2]$ such vectors. This is enough for us but it's also easy and well known how to use these $n/2$ orthonormal vectors to get $n$ orthonormal vectors with a somewhat worse lower bound on their norms.)

The map $T:\ell_2^n\to X$ defined by $Te_i=x_i$ is norm one. Note that
\[
1=\|T\|=\sup_{\|x^*\|_{X^*}\le 1} (\sum_{i=1}^n (x^*(x_i))^2)^{1/2}.
\]
When $X$ is isometric to a subspace of $\ell_\infty^N$ there are $N$ elements $x_j^*\in B_{X^*}$ such that, for all $x\in X$, $\|x\|=\max_{1\le j\le N}x_j^*(x)$. From this it is easy to deduce that
\[
\sup_{\|x^*\|_{X^*}\le 1} (\sum (x^*(x_i))^2)^{1/2}=\max_{1\le j\le N}(\sum_{i=1}^n (x_j^*(x_i))^2)^{1/2}.
\]
Denoting $a_i(j)=x_j^*(x_i)$ we get (\ref{eq:2summ}).

\begin{re}\label{re:ex}
Here we would like to sugget an approach toward showing that the dichotomy conjecture fails and maybe even that one can't get below the estimate $n>(\log N)^2$ in Theorem \ref{thm:main}.
\end{re}

Let $X$ and $Y$ be two $l$ dimensional normed spaces. Put $n=l^2$ and $N=36^l$. Let $\{x_i\}_{i=1}^{6^l}$ be a $1/2$ net in the sphere of $X$ and $\{y^*_i\}_{i=1}^{6^l}$ be a $1/2$ net in the sphere of $Y^*$. Note that for every $T:X\to Y$,
\[
\max_{1\le i,j\le 6^l}y^*_i(Tx_j)\le\|T\|\le 4\max_{1\le i,j\le 6^l}y^*_i(Tx_j).
\]
Consequently, $B(X,Y)$, the space of operators from $X$ to $Y$ with the operator norm, $4$-embeds into $\ell_\infty^N$. Note that $\dim(B(X,Y))=n\sim (\log N)^2$.

(Un)fortunately, $B(X,Y)$ cannot serve as a negative example since it always contains $\ell_\infty$-s with dimension going to infinity with $N$. This was pointed out to us by Bill Johnson. Indeed, by Dvoretzky's theorem, $\ell_2^k$ 2-embeds into $Y$ and into $X^*$, for some $k$ tending to infinity with $n$. Let $I$ denote the first embedding and $Q$ be the adjoint of the second embedding. It is then easy to see that $T\to ITQ$ is a 4-embedding of $B(\ell_2^k,\ell_2^k)$ into $B(X,Y)$. Finally, $B(\ell_2^k,\ell_2^k)$ contains isometrically $\ell_\infty^k$.

However, to get a negative answer to the dichotomy problem, it is enough to find $n$ dimensional $X$ and $Y$ and a subspace $Z$ of $B(X,Y)$ of dimension $m$ with $m/n$ tending to infinity with $n$ which has good cotype, i.e., if $Z$ contains a $2$-isomorph of $\ell_\infty^k$ then $k$ is bounded by a universal constant.
If one can find such an example with $m\ge cn^2$ for some universal positive constant $c$ then it will even show that one can't get below the estimate $n>(\log N)^2$ in Theorem \ref{thm:main}.

%\vfill\eject

 \begin{tabular}{l}
G.~Schechtman\\
Department of Mathematics\\
Weizmann Institute of Science\\
Rehovot, Israel\\
{\tt gideon@weizmann.ac.il}\\
\end{tabular}

\bigskip

\begin{tabular}{l}
N. Tomczak--Jaegermann\\
Department of Mathematics\\
University of Alberta\\
Edmonton, Canada\\
{\tt nicole.tomczak@ualberta.ca}\\
\end{tabular}


\begin{thebibliography}{99}

%\bibitem[B1]{b1}  Bourgain, Jean . New classes of Lp-spaces.
%Lecture Notes in Mathematics, 889. Springer-Verlag, Berlin-New York,  1981.

\bibitem[B1]{b2} Bourgain, J.,  Subspaces of $L_N^\infty$, arithmetical diameter and Sidon sets.
 Probability in Banach spaces, V (Medford, Mass., 1984),
 96--127, Lecture Notes in Math., 1153, Springer, Berlin,  1985.

\bibitem[B2]{b3} Bourgain, J.,  Bounded orthogonal systems and the $\Lambda(p)$-set problem.
 Acta Math.  162  (1989),  no. 3-4, 227--245.

% \bibitem[BT]{bt} Bourgain, J. ;  Tzafriri, L.  Restricted invertibility of matrices and %applications.
% Analysis at Urbana, Vol. II (Urbana, IL, 1986�1987),
% 61--107, London Math. Soc. Lecture Note Ser., 138, Cambridge Univ. Press, Cambridge,  1989.

\bibitem[FJ]{fj}   Figiel, T.;  Johnson, W. B.,  Large subspaces of $l_\infty^n$ and estimates of %the
 Gordon-Lewis constant.
 Israel J. Math.  37  (1980),  no. 1-2, 92--112.

\bibitem[GG]{gg} Garling, D. J. H.; Gordon, Y., Relations between some constants associated with finite dimensional Banach spaces. Israel J. Math. 9 (1971), 346--361.

\bibitem[MS]{ms} Milman, V. D.;  Schechtman, G., Asymptotic theory of finite-dimensional normed spaces.
With an appendix by M. Gromov.
Lecture Notes in Mathematics, 1200. Springer-Verlag, Berlin,  1986. {\rm viii}+156 pp.

\bibitem[P]{p}   Pisier, G.,  Remarques sur un r\'esultat non publi\'e de B. Maurey.
(French)  [[Remarks on an unpublished result of B. Maurey]]  Seminar on Functional Analysis, 1980--1981,
 Exp. No. V, 13 pp., \'Ecole Polytech., Palaiseau,  1981.


\bibitem[S]{s} Snobar, M. G., p-absolutely summing constants. (Russian) Teor. Funkci\u i Funkcional. Anal. i Prilo\v zen. No. 16 (1972), 38--41, 216.


\bibitem[T-J]{t-j} Tomczak-Jaegermann, N., Computing 2-summing norm with few vectors.
 Ark. Mat.  17  (1979),  no. 2, 273--277.

 \bibitem[T-J1]{t-j1} Tomczak-Jaegermann, N., Banach-Mazur distances and finite-dimensional operator ideals. Pitman Monographs and Surveys in Pure and Applied Mathematics, 38. Longman Scientific \& Technical, Harlow; copublished in the United States with John Wiley \& Sons, Inc., New York, 1989.


 \end{thebibliography}
\end{document}